\documentclass[11pt]{amsart}
\usepackage{amsmath,amsthm,latexsym,amssymb}
\usepackage{amsmath}
\usepackage{amsfonts}
\usepackage{amssymb}
\usepackage{graphicx}
\usepackage{color}
\providecommand{\U}[1]{\protect\rule{.1in}{.1in}}
\numberwithin{equation}{section}

\newtheorem{theorem}{Theorem}
\newtheorem{lemma}{Lemma}

\newtheorem{proposition}{Proposition}
\newtheorem{remark}{Remark}

\numberwithin{theorem}{section}
\numberwithin{corollary}{section}
\numberwithin{lemma}{section}
\numberwithin{definition}{section}
\numberwithin{proposition}{section}
\numberwithin{remark}{section}

\newcommand{\medint}{-\kern  -,375cm\int}

\begin{document}
\title[\bigskip]{Isoperimetric estimates for the first Neumann eigenvalue of Hermite
differential equations}
\author{F. Chiacchio$^{1}$ - G. Di Blasio$^{2}$}

\begin{abstract}
We provide isoperimetric Szeg\"{o}-Weinberger type inequalities for the first
nontrivial Neumann eigenvalue $\mu_{1}(\Omega)$ in Gauss space, where $\Omega$ is a possibly unbounded domain of $\mathbb{R}^{N}$.
Our main result consists in showing that among all sets of
$\mathbb{R}^{N}$ symmetric about the origin, having prescribed Gaussian measure, $\mu_{1}(\Omega)$ is
maximum if and only if $\Omega$ is the euclidean ball centered at the origin.

\vspace{1cm}


\vspace{0.3cm}

\textrm{\textsl{Mathematics Subject Classification:}}\noindent\ 35B45, 35P15, 35J70

\end{abstract}
\maketitle

\setcounter{footnote}{1} \footnotetext{Dipartimento di Matematica e
Applicazioni \textquotedblleft R. Caccioppoli\textquotedblright,
Universit\`{a} degli Studi di Napoli \textquotedblleft Federico
II\textquotedblright, Complesso Monte S. Angelo, via Cintia, 80126 Napoli,
Italy, e-mail: francesco.chiacchio@unina.it} \setcounter{footnote}{2}
\footnotetext{Dipartimento di Matematica, Seconda Universit\`{a} degli Studi
di Napoli, via Vivaldi, Caserta, Italy, e-mail: giuseppina.diblasio@unina2.it}

\section{Introduction}

Let us consider the classical eigenvalue problem for the free membrane
\begin{equation}
\left\{
\begin{array}
[c]{ccc}
-\Delta u=\mu u & \text{in} & \Omega\\
&  & \\
\dfrac{\partial u}{\partial\nu}=0 & \text{on} & \partial\Omega,
\end{array}
\right.  \label{P_I}%
\end{equation}
where $\Omega$ is a smooth connected subset of $\mathbb{R}^{N}$ and $\nu$ is
the outward normal to $\partial\Omega.$

In\ \cite{KS} Kornhauser and Stakgold conjectured that among all planar simply
connected domains, with fixed measure, $\mu_{1}(\Omega)$, the first nontrivial
eigenvalue of (\ref{P_I}), achieves its maximum value if and only if $\Omega$
is a disk.

This conjecture was proved by Szeg\"{o} in \cite{Sz}, by means of tools from
complex analysis, in particular he used the invariance of Dirichlet integrals
under conformal transplantation.

Soon after (see \cite{W}) Weinberger generalized this result to any bounded
smooth domain $\Omega$ of $\mathbb{R}^{N}.$ Weinberger obtained from the
eigenfunctions of the unit ball of $\mathbb{R}^{N}$, $B_{1}$, test functions
admissible in the variational characterization of $\mu_{1}(\Omega)$. His idea
was to extend radially such eigenfunctions in $\mathbb{R}^{N}$, just setting
their value constant outside $B_{1}$. Via the so-called \textquotedblleft
center of mass\textquotedblright\ arguments, he obtained $N$ different
functions having mean value zero on $\Omega$. At this point he is allowed to
use all these functions as trial functions for $\mu_{1}(\Omega)$ and the
result is finally achieved by symmetrization arguments.

This last method turned out to be rather flexible. We recall indeed that,
adapting Weinberger arguments, similar inequalities for spaces of constant
sectional curvature are derived. For instance in \cite{cha1} and \cite{AB3} it
is shown that if $\Omega$ is a domain of $\mathbb{S}^{N},$ contained in a
hemisphere, then
\[
\mu_{1}(\Omega)\leq\mu_{1}(\Omega^{\sharp}),
\]
where $\Omega^{\sharp}$ is the cap (i.e. the geodesic ball in $\mathbb{S}^{N}
$) having the same measure as $\Omega$.

On the other hand in \cite{La} it is proved that the first nonzero Neumann
eigenvalue is maximal for the equilateral triangle among all triangles of
given perimeter, and hence among all triangles of given area.

For further references see, e.g., the monographs \cite{Ba}, \cite{cha1},
\cite{Ke} and the survey paper \cite{A}.

The present paper deals with the Neumann eigenvalue problem in Gauss space. More precisely we study the problem
\begin{equation}
\left\{
\begin{array}
[c]{ll}
-\sum\limits_{i=1}^{N}\dfrac{\partial}{\partial x_{i}}\left(  \varphi_{N}\left(  x\right)
\dfrac{\partial u}{\partial x_{i}}\right)  =\mu\varphi_{N}\left(  x\right)
u & \text{ in }\Omega,\\
& \\
\dfrac{\partial u}{\partial\nu}=0 & \text{ on }\partial\Omega,
\end{array}
\right.  \label{P2_I}
\end{equation}
where $\Omega$ is a Lipschitz domain of $\mathbb{R}^{N}$, $N\geq1$ and
$\varphi_{N}\left(  x\right)  =(2\pi)^{-\frac{N}{2}}
e^{-\frac{|x|^{2}}{2}}$ is the density of normalized $N$-dimensional
Gaussian measure $d\gamma_{N}=\varphi_{N}\left(  x\right)  dx$.

Since the first half of the last century problems of the type (\ref{P2_I})
have attracted attention among both pure mathematicians and physicists. There
is indeed a tight connection between the eigenvalues of (\ref{P2_I}) and the
energy levels of the $N-$dimensional quantum harmonic oscillator. Related
references are the classical Courant-Hilbert monographs \cite{C-H} (see also
\cite{Fl}). On the other hand the interest in probability is motivated, for
instance, from the fact that the differential operator $L=-\Delta+x\cdot
\nabla$ appearing at the left hand side of (\ref{P2_I}) is the generator of
the Ornstein-Uhlenbeck semigroup, see, e.g., \cite{Bog} and the references
therein. Finally problems of the type (\ref{P2_I}) are related to some
functional inequalities as the well known Gross's Theorem on the Sobolev
Logarithmic embedding (see \cite{G}, \cite{Ad},\ \cite{PT}, \cite{CP} and
\cite{FP}).

If $\Omega$ is the whole space $\mathbb{R}^{N}$ the eigenfunctions of
(\ref{P2_I}) are the Hermite polynomials. If instead $\Omega\subsetneqq
\mathbb{R}^{N}$ then sharp estimates for the eigenvalues and eigenfunctions of
(\ref{P2_I}) with zero boundary conditions are contained, e.g., in \cite{Eh1},
\cite{BCM} and \cite{Be}.

Our aim is to prove isoperimetric Szeg\"{o}-Weinberger type inequalities for
the first eigenvalue of (\ref{P2_I}) that, with an abuse of notation, we still
denote by $\mu_{1}\left(  \Omega\right)  $.

In this setting it appears natural to maximize $\mu_{1}\left(  \Omega\right)
$ keeping fixed the Gaussian measure of $\Omega$. We recall that the Gaussian
measure in $\mathbb{R}^{N}$ can be obtained as a limit, as $k$ goes to
infinity, of normalized surface measures on $\mathbb{S}_{\sqrt{k}}^{k+N+1},$
the sphere in $\mathbb{R}^{k+N+2}$ of radius $\sqrt{k},$ (a process known in
literature as \textquotedblleft Poincar\'{e} limit\textquotedblright). Using
this limit process many properties for the Gauss space (i.e. $\mathbb{R}^{N}$
equipped with the measure $d\gamma_{N}$) can be deduced from analogous
properties which hold true for the sphere. One of the most remarkable example
is the Gaussian isoperimetric inequality, which asserts that among all
subsets $G$ of $\mathbb{R}^{N}$ with fixed Gaussian measure, the half-spaces
achieve the smallest Gaussian perimeter (see, e.g., \cite{Bo}, \cite{S} and
\cite{C-K}). We recall indeed, see, e.g., \cite{E}, that the half-spaces are
the \textquotedblleft Poincar\'{e} limit" of the caps, which are in turn the optimal sets in
the isoperimetric problem on the sphere. Another example is the Faber-Krahn
type inequality for Gauss space: the first Dirichlet Gaussian eigenvalue is
minimum on the half-space (see, e.g., \cite{E} and \cite{BCM}).

There are therefore two clues that might lead one to think that the half-space
would be a good candidate to maximize $\mu_{1}$. One reason is that the caps
maximize the first Neumann eigenvalue on the sphere, the other is that in all
the classical situations, described before, it is always the isoperimetric set
to maximize $\mu_{1}$.

This phenomenon here does not occur.

In one dimension we provide a detailed description of the behavior of $\mu
_{1}$. Let $\Omega=(a,b)$ with $-\infty\leq a<b\leq+\infty$ and $\gamma
_{1}(a,b)=L\in(0,1).$ We prove that $\mu_{1}(a,b)$ is minimum when the
interval reduces to an half-line, it is maximum when it is centered at the
origin and finally $\mu_{1}(a,b)$ is strictly monotone as $(a,b)$ slides
between these extreme positions. Therefore the set which gives the highest
eigenvalue is the one which maximizes the weighted perimeter and vice versa.

Our main result, which goes in the same direction of the previous one,
concerns the $N-$ dimensional case. We show that among all connected
 and possibly unbounded domain $\Omega$  of $\mathbb{R}^{N}$,
 symmetric about the origin and with fixed Gaussian measure, $\mu_{1}(\Omega)$ achieves its maximum value if
and only if $\Omega$ is the euclidean ball.

Since, obviously, the half-spaces are not symmetric about the origin, the above
result cannot exclude the possibility that the half-spaces maximize $\mu
_{1}\left(  \Omega\right)  $ in dimension greater then one. We are able to
exclude this possibility providing a suitable counterexample, see Remark
\ref{remark}.

We finally note that, as well known, a distinctive feature of Gauss measure is that its density
has both radial symmetry and product structure. Hence in the problem under consideration the former feature prevails on the latter.

\section{Notation and preliminary results}

Here and in the sequel $\Omega$ will denote a connected, smooth, open subset
of $\mathbb{R}^{N}$ such that $\gamma_{N}\left(  \Omega\right)  :=\int
_{\Omega}d\gamma_{N}<1.$ The natural functional space associated to problem
(\ref{P2_I}) is $H^{1}\left(  \Omega,\gamma_{N}\right)  $ which is the
weighted Sobolev space defined as follows%
\[
H^{1}\left(  \Omega,\gamma_{N}\right)  =\left\{  u\in W_{\text{loc}}
^{1,1}\left(  \Omega\right)  :\left(  u,\left\vert Du\right\vert \right)  \in
L^{2}\left(  \Omega,\gamma_{N}\right)  \times L^{2}\left(  \Omega,\gamma
_{N}\right)  \right\}  ,
\]
endowed with the norm
\begin{equation}
\left\Vert u\right\Vert _{H^{1}\left(  \Omega,\gamma_{N}\right)  }=\left\Vert
u\right\Vert _{L^{2}\left(  \Omega,\gamma_{N}\right)  }+\left\Vert
Du\right\Vert _{L^{2}\left(  \Omega,\gamma_{N}\right)  }=
\left(\int_{\Omega}u^{2}d\gamma_{N}\right)^{\frac{1}{2}}+
\left(\int_{\Omega}\left\vert Du\right\vert ^{2}d\gamma_{N}\right)^{\frac{1}{2}}.
\label{norm_H^1}
\end{equation}
In \cite{FP}, among other things, it is proved that the subspace of
$H^{1}\left(  \Omega,\gamma_{N}\right)  $ made of those functions having mean
value zero in $\Omega$ it is compactly embedded in $L^{2}\left(  \Omega
,\gamma_{N}\right)  .\,$\ This circumstance allows us to use standard spectral
theory for self-adjoint compact operator. In particular the variational
characterization of $\mu_{1}\left(  \Omega\right)  $ will be used throughout
\begin{equation}
\mu_{1}\left(  \Omega\right)  =\min_{\substack{v\neq0\\\int_{\Omega}
vd\gamma_{N}=0}}\frac{\int_{\Omega}\left\vert Dv\right\vert ^{2}d\gamma_{N}
}{\int_{\Omega}v^{2}d\gamma_{N}}. \label{def_auto}
\end{equation}
We recall, see, e.g., \cite{C-H}, that when $\Omega=\mathbb{R}^{N}$ the
eigenfunctions to problem (\ref{P2_I}) are combinations of homogeneous Hermite
polynomials. The Hermite polynomials in one variable are defined by
\begin{equation}
H_{n}(t)=(-1)^{n}e^{t^{2}/2}\frac{d^{n}}{dt^{n}}e^{-t^{2}/2},\text{ \ \ }
n\in\mathbb{N}\cup\left\{  0\right\}  , \label{Hermite}
\end{equation}
and they constitute a complete set of eigenfunctions to problem (\ref{P2_I})
with $\Omega=\mathbb{R},$ more precisely it holds
\[
-(\varphi_{1}\left(  t\right)  H_{n}^{\prime}(t))^{\prime}=n\varphi_{1}\left(
t\right)  H_{n}(t).
\]
Since $\Omega$ is a smooth set,$\,$its Gaussian perimeter is simply given by
\[
P_{\gamma_{N}}\left( \Omega \right)  =\int_{\partial \Omega}\gamma_{N}\left(  x\right)
\ \mathcal{H}_{N-1}\left(  dx\right)  ,
\]
where $\mathcal{H}_{N-1}\left(  x\right)  $ is the $(N-1)-$dimensional
Hausdorff measure.

As already mentioned in the introduction, for the Gaussian measure an
isoperimetric inequality holds true. Consider the half-space
\begin{equation}
\Omega^{\bigstar}=\left\{  x\in\mathbb{R}^{N}:x_{1}>\Phi^{-1}\left(
\gamma_{N}\left(  \Omega\right)  \right)  \right\}  , \label{Omega_star}
\end{equation}
where $\Phi\left(  t\right)  $ is the complementary error function
\begin{equation}
\Phi\left(  t\right)  =\frac{1}{\sqrt{2\pi}}\int_{t}^{\infty}e^{-\frac{s^{2}
}{2}}ds. \label{phi}
\end{equation}
In other words $\Omega^{\bigstar}$ is the half-space orthogonal to the $x_{1}
$-axis having the same Gaussian measure as $\Omega$.

The isoperimetric inequality for Gaussian measure (see \cite{S}, \cite{Bo},
\cite{Eh1} and \cite{C-K}) states that
\begin{equation}
P_{\gamma_{N}}\left(  \Omega\right)  \geq P_{\gamma_{N}}\left(  \Omega
^{\bigstar}\right)  , \label{dis isop}
\end{equation}
where equality holds in (\ref{dis isop}) if and only if $\Omega=\Omega
^{\bigstar},$ modulo a rotation.

Now we recall a few definitions and properties about Gaussian rearrangement,
whose notion was introduced by Ehrhard in \cite{Eh1}. For exhaustive treatment
on rearrangements we refer, e.g., to \cite{Ba}, \cite{CR}, \cite{Ka} and
\cite{Rs}.

Let $u:x\in\Omega\rightarrow\mathbb{R}$ be a measurable function. We denote by
$\mu(t)$ the distribution function of $\left\vert u(x)\right\vert $ i.e.
\[
\mu(t)=\gamma_{N}\left(  \left\{  x\in\Omega:\left\vert u(x)\right\vert
>t\right\}  \right)  ,\quad t\geq0,
\]
while the decreasing rearrangement and the increasing rearrangement of $u$,
with respect to the Gaussian measure, are defined respectively by
\[
u^{\ast}\left(  s\right)  =\inf\left\{  t\geq0:\mu\left(  t\right)  \leq
s\right\}  \text{,}\quad s\in\left]  0,\gamma_{N}(\Omega)\right]
\]
and%
\[
u_{\ast}\left(  s\right)  =u^{\ast}\left(  \gamma_{N}\left(  \Omega\right)
-s\right)  \text{,}\quad s\in\left[  0,\gamma_{N}(\Omega)\right[  .
\]
Finally $u^{\bigstar}$, the Gaussian rearrangement of $u$, is given by%
\[
u^{\bigstar}\left(  x\right)  =u^{\ast}\left(  \Phi\left(  x_{1}\right)
\right)  ,\quad x\in\Omega^{\bigstar}.
\]
By its very definition $u^{\bigstar}$ depends on one variable only and it is
an increasing function, therefore its level sets are parallel half-spaces.
Since, by definition, $u$ and $u^{\bigstar}$ are equimisurable, Cavalieri's
principle ensures%
\[
\left\Vert u\right\Vert _{L^{p}\left(  \Omega,\gamma_{N}\right)  }=\left\Vert
u^{\bigstar}\right\Vert _{L^{p}\left(  \Omega^{\bigstar},\gamma_{N}\right)
},\text{ }\forall p\geq1.
\]
We will also make use of the Hardy-Littlewood inequality, which states that
\begin{equation}
\int_{0}^{\gamma_{N}\left(  \Omega\right)  }u^{\ast}\left(  s\right)  v_{\ast
}\left(  s\right)  ds\leq\int_{\Omega}\left\vert u\left(  x\right)  v\left(
x\right)  \right\vert d\gamma_{N}\leq\int_{0}^{\gamma_{N}\left(
\Omega\right)  }u^{\ast}\left(  s\right)  v^{\ast}\left(  s\right)  ds.
\label{HL}
\end{equation}
We finally recall the Polya-Szeg\"{o} principle which asserts that the
weighted $L^{2}$ norms of a nonnegative function vanishing on $\partial\Omega$
decreases under Gaussian symmetrization. More precisely let $H_{0}^{1}\left(
\Omega,\gamma_{N}\right)  $ be the closure of $C_{0}^{\infty}\left(
\Omega\right)  $ in $H^{1}\left(  \Omega,\gamma_{N}\right)  .$ It holds that
\[
\int_{\Omega}\left\vert Du\left(  x\right)  \right\vert ^{2}d\gamma_{N}
\geq\int_{\Omega^{\bigstar}}\left\vert Du^{\bigstar}\left(  x\right)
\right\vert ^{2}d\gamma_{N},
\]
for any nonnegative $u$ in $H_{0}^{1}\left(  \Omega,\gamma_{N}\right)  $.

\section{The one-dimensional case}

Let $a,b\in\overline{\mathbb{R}}$ with $-\infty\leq a<b\leq+\infty$ and
$\gamma_{1}(a,b)<1.$ In this case problem (\ref{P2_I}) becomes%
\begin{equation}
\left\{
\begin{array}
[c]{ll}
-u^{\prime\prime}+xu^{\prime}=\mu u & \text{ in }(a,b)\\
& \\
u^{\prime}(a)=u^{\prime}(b)=0. &
\end{array}
\right.  \label{1d_N}
\end{equation}
We will denote by $\mu_{1}(a,b)$ the first nontrivial eigenvalue of
(\ref{1d_N}), clearly its value is given by
\begin{equation}
\mu_{1}(a,b)=\min_{u\neq0:\int_{a}^{b}ud\gamma_{1}=0}\frac{\int_{a}%
^{b}(u^{\prime})^{2}d\gamma_{1}}{\int_{a}^{b}u^{2}d\gamma_{1}}.
\label{lamda_1(a,b)}
\end{equation}
Here we are interested in studying the behavior of $\mu_{1}(a,b)$ when the
interval $(a,b)$ slides along the x-axis, keeping fixed its Gaussian measure.
In other words, we impose the constraint
\begin{equation}
\gamma_{1}(a,b)=L\in(0,1). \label{constraint_L}%
\end{equation}
Obviously, under these conditions, $b$ can be expressed in terms of $a$ as
follows
\begin{equation}
b(a)=\sqrt{2}\operatorname{erf}^{-1}\left[  2L+\operatorname{erf}\left(
\frac{\sqrt{2}}{2}a\right)  \right]  , \label{b(a)}
\end{equation}
where
\[
\operatorname{erf}(x)=\frac{2}{\sqrt{\pi}}\int_{0}^{x}e^{-t^{2}}dt
\]
is the error function.

Since condition (\ref{constraint_L}) is in force, the function
\begin{equation}
f:a\in\mathbb{R}\rightarrow\mu_{1}(a,b(a)) \label{f}
\end{equation}
is defined on the interval $\left[  -\infty,\sqrt{2}\operatorname{erf}
^{-1}(1-2L)\right]  $ and it is even with respect to $x=-\sqrt{2}
\operatorname{erf}^{-1}(L)$.

The following result holds.

\begin{theorem}
\label{Th-OneDim}Let $L\in(0,1)$ and let $-\infty\leq a<b\leq+\infty,$ with
$\gamma_{1}(a,b\left(  a\right)  )=$ $L$. Then
\begin{equation}
\min_{a}\mu_{1}(a,b(a))=\mu_{1}(-\infty,\sqrt{2}\operatorname{erf}%
^{-1}(2L-1))=\mu_{1}(\sqrt{2}\operatorname{erf}^{-1}(1-2L),+\infty),
\label{min}
\end{equation}
and
\begin{equation}
\max_{a}\mu_{1}(a,b(a))=\mu_{1}(-\sqrt{2}\operatorname{erf}^{-1}(L),\sqrt
{2}\operatorname{erf}^{-1}(L)). \label{max}
\end{equation}
Furthermore the function $f$ defined in (\ref{f}) is increasing in the
interval $\left[  -\infty,-\sqrt{2}\operatorname{erf}^{-1}(L)\right]  $.
\end{theorem}

\textbf{Proof. }We denote by $\lambda_{1}(a,b(a))$ the first eigenvalue of the
problem%
\begin{equation}
\left\{
\begin{array}
[c]{ll}%
-v^{\prime\prime}+xv^{\prime}=\lambda v & \text{ in }(a,b(a))\\
& \\
v(a)=v(b(a))=0. &
\end{array}
\right.  \label{1d_D}
\end{equation}
It is easy to verify that
\begin{equation}
\lambda_{1}(a,b(a))=\mu_{1}(a,b(a))-1. \label{D=N-1}
\end{equation}

Indeed let $u_{1}$ be an eigenfunction corresponding to $\mu _{1}(a,b(a)),$
then $v =u_{1}^{\prime }$ satisfies (\ref{1d_D}) with $\lambda =$ $\mu
_{1}(a,b(a))-1.$ This means that $\mu _{1}(a,b(a))\geq \lambda _{1}(a,b(a))+1
$. It remains to prove the converse inequality. To this aim consider the
function $\psi (x,x_{0})=\int_{x_{0}}^{x}v_{1}(\sigma )d\sigma $, where $v_{1}$
is an eigenfunction of (\ref{1d_D}) corresponding to $\lambda _{1}(a,b(a))$
and $x_{0}\in (a,b)$ is chosen such that $\int_{a}^{b}\psi (x,x_{0})d\gamma
_{1}=0$. The function $\psi $ verifies
\begin{equation*}
\frac{\int_{a}^{b}\left( \psi ^{\prime }\right) ^{2}d\gamma _{1}}{
\int_{a}^{b}\left( \psi \right) ^{2}d\gamma _{1}}=\lambda _{1}(a,b(a))+1,
\end{equation*}
and $\psi ^{\prime }(a)=\psi ^{\prime }(b)=0$. Therefore $\mu
_{1}(a,b(a))\leq \lambda _{1}(a,b(a))+1$, which is the claim.

Since they differ by a constant, in place of the Neumann eigenvalue we can
equivalently study the behavior of the Dirichlet eigenvalue.

As a first consequence of this observation we note that the Faber-Krahn
inequality for Gaussian measure (see \cite{E} and \cite{BCM}) directly gives
(\ref{min}).

The isoperimetric properties of the half-space (see, e.g., \cite{Bo} and
\cite{S}) reads as follows%
\[
\min_{a}P_{\gamma_{1}}(a,b(a))=P_{\gamma_{1}}(-\infty,\sqrt{2}
\operatorname{erf}^{-1}(2L-1))=P_{\gamma_{1}}(\sqrt{2}\operatorname{erf}
^{-1}(1-2L),+\infty).
\]
A straightforward application of Lagrange multipliers rule tells us that the
function $P_{\gamma_{1}}(a,b)$ admits just one stationary point on the
constraint $\gamma_{1}(a,b)-L=0$. Moreover, as it is immediate to verify, such
a point occurs at $a=-b=-\sqrt{2}$ $\operatorname{erf}^{-1}(L).$ Now since the
function $P_{\gamma_{1}}(a,b(a))$ is smooth on the interval $(-\infty,\sqrt
{2}\operatorname{erf}^{-1}(1-2L)),$ from (\ref{min}) we get%
\begin{equation}
\max_{a}P_{\gamma_{1}}(a,b(a))=P_{\gamma_{1}}(-\sqrt{2}\operatorname{erf}%
^{-1}(L),\sqrt{2}\operatorname{erf}^{-1}(L)). \label{max_per}%
\end{equation}
These considerations allow us to say that
\begin{equation}
\frac{d}{da}P_{\gamma_{1}}(a,b(a))>0\text{ \ \ \ \ }\forall a\in
(-\infty,-\sqrt{2}\operatorname{erf}^{-1}(L)) \label{P_incr}%
\end{equation}
and by symmetry reasons%
\begin{equation}
\frac{d}{da}P_{\gamma_{1}}(a,b(a))<0\text{ \ \ \ \ }\forall a\in(-\sqrt
{2}\operatorname{erf}^{-1}(L),\sqrt{2}\operatorname{erf}^{-1}(1-2L)).
\label{P_dec}%
\end{equation}
Now we can finally turn our attention on the monotonicity properties of the
eigenvalue $\mu_{1}(a,b(a)).$ Let $a_{1},a_{2}\in(-\infty,-\sqrt
{2}\operatorname{erf}^{-1}($ $L))$ with $a_{1}<a_{2}.$ Our aim is to prove
that%
\begin{equation}
\mu_{1}(a_{1},b(a_{1}))<\mu_{1}(a_{2},b(a_{2})) \label{auto_inc}%
\end{equation}
or equivalently%
\[
\lambda_{1}(a_{1},b(a_{1}))<\lambda_{1}(a_{2},b(a_{2})).
\]
Let us denote by $\phi_{i}(x)$, with $i=1,2$, the first Dirichlet
eigenfunctions corresponding to $\lambda_{1}(I_{i})$, where $I_{i}%
=(a_{i},b(a_{i}))$, with $i=1,2$, normalized in such a way that they are
positive and
\[
\int_{a_{i}}^{b(a_{i})}\phi_{i}^{2}(x)d\gamma_{1}=1.
\]

For any fixed $t\in(0,1)$ we denote by $I_{2}^{t}$ the set $\left\{  x\in
I_{2}:\phi_{2}(x)>t\right\}  .$ From the level sets of $\phi_{2}(x)$ we want
to build a function defined in $I_{1}$ admissible as test function for
$\lambda_{1}(I_{1})$. This auxiliary function, denoted with $\widetilde{\phi
}(x)$, is the function uniquely defined by the following relationships

\begin{itemize}
\item[(i)] $\widetilde{\phi}:x\in I_{1}\rightarrow\left[  0,\max\phi
_{2}\right]  ,$

\item[(ii)] $\gamma_{1}\{x:\widetilde{\phi}(x)>t\}=\gamma_{1}\left\{
x:\phi_{2}(x)>t\right\}  ,$ $\ \ \ \ \ \forall t\in\left[  0,\max\phi
_{2}\right]  ,$

\item[(iii)] $\{x:\widetilde{\phi}(x)>t\}$ are intervals $\left(
\widetilde{a}_{t},\widetilde{b}_{t}\right)  $, denoted with $\overset{\sim
}{I_{t}},$ centered at $\frac{a_{1}+b(a_{1})}{2},$ $\ \ \forall t\in\left[
0,\max\phi_{2}\right]  .$
\end{itemize}

By construction $\widetilde{\phi}$ is even with respect to $\frac
{a_{1}+b(a_{1})}{2}$ and it is increasing in $\left(  a_{1},\frac
{a_{1}+b(a_{1})}{2}\right)  .$ Furthermore it is equimeasurable with $\phi
_{2}$, therefore $\mu_{\widetilde{\phi}}(t)=\mu_{\phi_{2}}(t)$ and%
\[
\int_{a_{1}}^{b(a_{1})}\widetilde{\phi}^{2}d\gamma_{1}=\int_{a_{2}}^{b(a_{2}%
)}\phi_{2}^{2}d\gamma_{1}=1.
\]
Coarea formula and Cauchy--Schwarz inequality\ ensure that%
\begin{align}
\lambda_{1}(I_{2})  &  =\frac{1}{\sqrt{2\pi}}\int_{a_{2}}^{b(a_{2})}\left(
\frac{d\phi_{2}}{dx}\right)  ^{2}e^{-\frac{x^{2}}{2}}dx=\frac{1}{\sqrt{2\pi}%
}\int_{0}^{\max\phi_{2}}\left(  \int_{\{\phi_{2}=t\}}\left\vert \frac
{d\phi_{2}}{dx}\right\vert e^{-\frac{x^{2}}{2}}d\mathcal{H}^{0}\right)
dt\label{L1_D_1ineq}\\
&  \geq\frac{1}{\sqrt{2\pi}}\int_{0}^{\max\phi_{2}}\frac{\left(  \int
_{\{\phi_{2}=t\}}e^{-\frac{x^{2}}{2}}d\mathcal{H}^{0}\right)  ^{2}}%
{\int_{\{\phi_{2}=t\}}\left\vert \frac{d\phi_{2}}{dx}\right\vert
^{-1}e^{-\frac{x^{2}}{2}}d\mathcal{H}^{0}}dt=\int_{0}^{\max\phi_{2}}%
\frac{\left(  P_{\gamma_{1}}\left\{  \phi_{2}>t\right\}  \right)  ^{2}}%
{-\mu_{\phi_{2}}^{\prime}(t)}dt.\nonumber
\end{align}
At this point we note that by (\ref{P_incr}) and by the construction of
$\widetilde{\phi}$ we have
\begin{equation}
P_{\gamma_{1}}\left\{  \phi_{2}>t\right\}  >P_{\gamma_{1}}\{\widetilde{\phi
}>t\}=P_{\gamma_{1}}\left(  \overset{\sim}{I_{t}}\right)  ,\text{
\ \ \ \ }\forall t\in(0,\max\phi_{2}) \label{P>P}%
\end{equation}
and%
\begin{equation}
\mu_{\phi_{2}}^{\prime}(t)=\mu_{\widetilde{\phi}}^{\prime}(t),\text{
\ \ \ \ a.e. }t\in(0,L). \label{mu'=mu'}%
\end{equation}

So by (\ref{L1_D_1ineq}), (\ref{P>P}) and (\ref{mu'=mu'}) we have%
\begin{equation}
\lambda_{1}(I_{2})>\int_{0}^{\max\widetilde{\phi}}\frac{\left(  P_{\gamma_{1}%
}\{\widetilde{\phi}>t\}\right)  ^{2}}{-\mu_{\widetilde{\phi}}^{\prime}(t)}dt.
\label{L1_D_2ineq}%
\end{equation}
Since the function $\widetilde{\phi}$ is, by construction, even with respect
to $\frac{a_{1}+b(a_{1})}{2}$ we have%
\[
\left\vert \frac{d\widetilde{\phi}}{dx}(\widetilde{a}^{t})\right\vert
=\left\vert \frac{d\widetilde{\phi}}{dx}(\widetilde{b}^{t})\right\vert ,\text{
\ \ \ \ a.e. }t\in(0,\max\phi_{2}),
\]
and therefore the Cauchy--Schwarz inequality used in (\ref{L1_D_1ineq}) for
$\widetilde{\phi}$ reduces to an equality. This consideration together with
(\ref{L1_D_2ineq}),\ yields%
\begin{equation}
\lambda_{1}(I_{2})\geq\frac{1}{\sqrt{2\pi}}\int_{a_{1}}^{b(a_{1})}\left(
\frac{d\widetilde{\phi}}{dx}\right)  ^{2}e^{-\frac{x^{2}}{2}}dx\geq\lambda
_{1}(I_{1}). \label{L1D>L1D}%
\end{equation}
That is the claim (\ref{auto_inc}). Note finally that if $\lambda_{1}%
(a_{1},b(a_{1}))=$ $\lambda_{1}(a_{2},b(a_{2}))$ then all the above
inequalities reduce to equalities. In particular equality in (\ref{P>P}) implies
that $a_{1}=a_{2\text{ }}$ and $b\left(  a_{1}\right)  =b\left(  a_{2}\right)
.$

\ \ \ \ \ \ \ \ \ \ \ \ \ \ \ \ \ \ \ \ \ \ \ \ \ \ \ \ \ \ \ \ \ \ \ \ \ \ \ \ \ \ \ \ \ \ \ \ \ \ \ \ \ \ \ \ \ \ \ \ \ \ \ \ \ \ \ \ \ \ \ \ \ \ \ \ \ \ \ \ \ \ \ \ \ \ \ \ \ \ \ \ \ \ \ \ \ \ \ \ \ \ \ \ \ \ \ \ \ \ \ \ \ \ \ \ \ \ \ \ \ \ \ \ \ \ \ \ \ \ \ \ \ \ \ \ \ \ \ \ \ $\square
$

\begin{remark}
\label{Deriv_atuo}Theorem \ref{Th-OneDim}, together with the shape derivative
formula for one-dimensional Neumann eigenvalues, allows to get some
qualitative information on $u_{1}$. Let us consider two smooth functions
$a(t)$ and $b(t),$ such that $\gamma_{1}\left(  a\left(  t\right)  ,b\left(
t\right)  \right)  =L$ and $(\left(  a\left(  0\right)  ,b\left(  0\right)
\right)  =(a,b).$ Let us denote by $\mu_{1}\left(  t\right)  =\mu_{1}\left(
a\left(  t\right)  ,b\left(  t\right)  \right)  $ the first eigenvalue of
problem
\[
\left\{
\begin{array}
[c]{ll}
-\dfrac{d^{2}}{dx^{2}}u\left(  x,t\right)  +x\dfrac{d}{dx}u\left(  x,t\right)
=\mu\left(  t\right)  u\left(  x,t\right)  & \text{ in }(a\left(  t\right)
,b\left(  t\right)  )\\
\left.  \dfrac{d}{dx}u\left(  x,t\right)  \right\vert _{x=a\left(  t\right)
,b\left(  t\right)  }=0, &
\end{array}
\right.
\]
and by $u_{1}\left(  x,t\right)  $ a corresponding eigenfunction such that
$\int_{a(t)}^{b(t)}u_{1}^{2}\left(  x,t\right)  d\gamma_{1}=1.$ Then, see,
e.g., \cite{H1, H2}, it is easy to verify that
\begin{equation}
\mu_{1}^{\prime}\left(  0\right)  =\mu_{1}(a,b)e^{-\frac{a^{2}}{2}}\left(
u^{2}\left(  a\right)  -u^{2}\left(  b\right)  \right) .
\label{Derivata_auto}
\end{equation}

Therefore if $a<-\sqrt{2}\operatorname{erf}^{-1}($ $L)$ then, by Theorem
\ref{Th-OneDim}, we have that $\left\vert u\left(  a\right)  \right\vert
>\left\vert u\left(  b\right)  \right\vert ,$ conversely if $a\in(-\sqrt
{2}\operatorname{erf}^{-1}(L),\sqrt{2}\operatorname{erf}^{-1}(1-2L))$ then
$\left\vert u\left(  a\right)  \right\vert <\left\vert u\left(  b\right)
\right\vert .$
\end{remark}

\section{The $N-$dimensional case}

Let us examine, by means of the separation of variables method, problem
(\ref{P2_I})\ when $\Omega$ is the ball of $\mathbb{R}^{N}$ centered at the
origin of radius $R$, throughout denoted by $B_{R}$, that is%
\begin{equation}
\left\{
\begin{array}
[c]{lll}%
-\Delta u+x\cdot Du=\mu u & \text{in} & B_{R}\\
&  & \\
\dfrac{\partial u}{\partial r}=0 & \text{on} & \partial B_{R}.
\end{array}
\right.  \label{P_Ball}%
\end{equation}
The equation in (\ref{P_Ball}) can be rewritten, using polar coordinates, as
\begin{equation}
\frac{1}{r^{N-1}}\frac{\partial}{\partial r}\left(  r^{N-1}\frac{\partial
u}{\partial r}\right)  +\frac{1}{r^{2}}\Delta_{\mathbb{S}^{N-1}}\left(
u|\mathbb{S}_{r}^{N-1}\right)  -r\frac{\partial u}{\partial r}+\mu
u=0,\label{laplaciano polare}%
\end{equation}
where $\mathbb{S}_{r}^{N-1}$ is the sphere of radius $r$ in $\mathbb{R}^{N}%
,$\textsl{ }$u|\mathbb{S}_{r}^{N-1}$\textsl{ }is the restriction of $u$ on
$\mathbb{S}_{r}^{N-1}$ and finally $\Delta_{\mathbb{S}^{N-1}}\left(
u|\mathbb{S}_{r}^{N-1}\right)  $ is the standard Laplace-Beltrami operator
relative to the manifold $\mathbb{S}_{r}^{N-1}.$

Setting $u\left(  x\right)  =Y\left(  \theta\right)  f\left(  r\right)  $\ in
equation (\ref{laplaciano polare}), where $\theta$ belongs to $\mathbb{S}%
_{1}^{N-1},$ we have%
\[
Y\frac{1}{r^{N-1}}\left(  r^{N-1}f^{\prime}\right)  ^{\prime}+\Delta
_{\mathbb{S}^{N-1}}Y\frac{f}{r^{2}}-Yrf^{\prime}+\mu Yf=0,
\]
and hence%
\begin{equation}
\frac{1}{r^{N-3}f}\left(  r^{N-1}f^{\prime}\right)  ^{\prime}-r^{3}
\frac{f^{\prime}}{f}+\mu r^{2}=-\frac{\Delta_{\mathbb{S}^{N-1}}Y}{Y}
=\overset{\_}{k}. \label{eq separate}
\end{equation}
As well known, see, e.g., \cite{Mu} and \cite{Cha}, the last equality is fulfilled if
and only if
\[
\overset{\_}{k}=k\left(  k+N-2\right)  \text{ \ \ with }k=\mathbb{N}
\cup\left\{  0\right\}  .
\]
Multiplying the left hand side of equation (\ref{eq separate}) by $\dfrac
{f}{r^{2}},$ we get
\[
f^{\prime\prime}+f^{\prime}\left(  \frac{N-1}{r}-r\right)  +\mu f-k\left(
k+N-2\right)  \dfrac{f}{r^{2}}=0\text{ \ in \ \ }\left(  0,R\right)  .
\]
The eigenfunctions are either purely radial%
\begin{equation}
u_{i}\left(  r\right)  =f_{0}\left(  \mu_{i};r\right)  ,\text{ if }k=0,
\label{Rad_eigen}
\end{equation}
or in the form

\begin{equation}
u_{i}\left(  r,\theta\right)  =f_{k}\left(  \mu_{i};r\right)  Y\left(
\theta\right)  ,\ \text{if}\ k\in\mathbb{N}. \label{Ang_eigen}
\end{equation}
The functions $f_{k},$ with $k\in\mathbb{N}\cup\left\{  0\right\}  ,$
\ clearly satisfy
\begin{equation}
\left\{  \left.
\begin{array}
[c]{l}%
f_{k}^{\prime\prime}+f_{k}^{\prime}\left(  \dfrac{N-1}{r}-r\right)  +\mu
_{i}f_{k}-k\left(  k+N-2\right)  \dfrac{f_{k}}{r^{2}}=0\text{ \ \ in
\ \ \ }\left(  0,R\right) \\
\\
f_{k}\left(  0\right)  =0,\text{ \ }f_{k}^{\prime}\left(  R\right)  =0.
\end{array}
\right.  \right.  \label{eq_for_f_k}
\end{equation}

In the sequel we will denote by $\tau_{n}(R),$ with $n\in\mathbb{N}
\cup\left\{  0\right\}  $, the sequence of eigenvalues of (\ref{P_Ball}) whose
corresponding eigenfunctions are purely radial, i.e. in the form
(\ref{Rad_eigen}) or equivalently solutions to problem (\ref{eq_for_f_k}) with
$k=0$. Clearly in this case the first eigenfunction is constant and the
corresponding eigenvalue $\tau_{0}(R)$ is trivially zero. We will denote by
$\nu_{n}(R)$, with $n\in\mathbb{N}$, the remaining eigenvalues of
(\ref{P_Ball}).

\begin{lemma}
\label{Lemma}It holds that
\begin{equation}
\nu_{1}(R)<\tau_{1}(R),\text{ \quad}\forall R>0. \label{vu1<tau2}
\end{equation}

\end{lemma}

\textbf{Proof. }We recall that $\tau_{1}=\tau_{1}(R)$ is the first nontrivial
eigenvalue of%
\begin{equation}
\left\{
\begin{array}
[c]{ll}%
g^{\prime\prime}+g^{\prime}\left(  \dfrac{N-1}{r}-r\right)  +\tau g=0 &
\mbox{ in }(0,R)\\
& \\
g^{\prime}\left(  0\right)  =g^{\prime}\left(  R\right)  =0, &
\end{array}
\right.  \label{prob_2}%
\end{equation}
and $\nu_{1}=\nu_{1}(R)$ is the first eigenvalue of
\begin{equation}
\left\{
\begin{array}
[c]{ll}%
w^{\prime\prime}+w^{\prime}\left(  \dfrac{N-1}{r}-r\right)  +\nu w-\left(
N-1\right)  \dfrac{w}{r^{2}}=0 & \mbox{ in }(0,R)\\
& \\
w\left(  0\right)  =w^{\prime}\left(  R\right)  =0. &
\end{array}
\right.  \label{prob_1}
\end{equation}
First of all\textbf{\ }we observe that the first eigenfunction $w_{1}$ of
(\ref{prob_1}) does not change its sign in $\left(  0,R\right)  $, thus we can
assume that $w_{1}>0$ in $\left(  0,R\right)  .$

Moreover $w_{1}^{\prime}\geq0$ in $\left(  0,R\right)  .$ Indeed, assume, by
contradiction, that we can find two values $r_{1}$, $r_{2},$ with $r_{1}<r_{2},$ such
that $w_{1}^{\prime\prime}\left(  r_{1}\right)  \leq0,$ $w_{1}^{\prime}\left(
r_{1}\right)  =0$ and $w_{1}^{\prime\prime}\left(  r_{2}\right)  \geq0,$
$w_{1}^{\prime}\left(  r_{2}\right)  =0.$ By evaluating the equation in
(\ref{prob_1})
\[
\frac{w_{1}^{\prime\prime}}{w_{1}}+\frac{w_{1}^{\prime}}{w_{1}}\left(
\frac{N-1}{r}-r\right)  +\nu_{1}-\frac{N-1}{r^{2}}=0
\]
at $r_{1}$ and $r_{2},$ we get
\[
\nu_{1}-\frac{N-1}{r_{2}^{2}}\leq0\text{ \ and \ \ }\nu_{1}-\frac{N-1}
{r_{1}^{2}}\geq0,
\]
that means $r_{1}\geq r_{2}$ and this is a contradiction.

On the other hand, the first nontrivial eigenfunction of problem
(\ref{prob_2}), $g_{1}=g_{1}(r)$, has mean value zero i.e.
\[
\int_{B_{R}}g_{1}d\gamma_{N}=\frac{N\omega_{N}}{\left(  2\pi\right)  ^{N/2}
}\int_{0}^{R}g_{1}(r)e^{-\frac{r^{2}}{2}}r^{N-1}dr=0,
\]
where, here and in the sequel, $\omega_{N}$ will denote the volume of the unit
ball in $\mathbb{R}^{N}$.

This implies that $g_{1}(r)$ must change its sign in $\left(  0,R\right)  $. Let
us suppose $g_{1}(r)>0$ in $\left(  0,r_{0}\right)  $ and $g_{1}\left(
r_{0}\right)  =0.$ We observe that $g_{1}^{\prime}(r)<0$ in $\left(
0,R\right)  .$ Moreover evaluating the equation of problem (\ref{prob_2}) at
$r_{0},$ we have
\begin{equation}
g_{1}^{\prime\prime}\left(  r_{0}\right)  +g_{1}^{\prime}\left(  r_{0}\right)
\left(  \frac{N-1}{r_{0}}-r_{0}\right)  =0. \label{eq_for_g''}%
\end{equation}
Now we consider the following intervals $J_{1}:=\left(  0,\sqrt{N-1}\right]
,$ $\ \ \ J_{2}:=\left(  \sqrt{N-1},\sqrt{N-1}+\frac{\pi}{\sqrt{8}}\right]  $
and $J_{3}:=\left(  \sqrt{N-1}+\frac{\pi}{\sqrt{8}},+\infty\right)  .$ Clearly
$\underset{i=1}{\overset{3}{\cup}}J_{i}=(0,+\infty)$ for any $N\in\mathbb{N}.$
The proof of (\ref{vu1<tau2}) requires different arguments depending on the
interval $J_{i}$ in which the radius $R$ of the ball $B_{R}$ lies.

\bigskip

\textbf{Case 1}: $R\in J_{1}=\left(  0,\sqrt{N-1}\right]  .$

Since $r_{0}<R\leq\sqrt{N-1},$ from (\ref{eq_for_g''}) we get
\begin{equation}
g_{1}^{\prime\prime}\left(  r_{0}\right)  \geq0. \label{g''(r_0)>0}%
\end{equation}
Moreover if we set $\psi=g_{1}^{\prime},$ then problem (\ref{prob_2}) becomes%
\[
\left\{
\begin{array}
[c]{ll}
\psi^{\prime\prime}+\psi^{\prime}\left(  \dfrac{N-1}{r}-r\right)  +\psi\left(
-\dfrac{N-1}{r^{2}}-1\right)  +\tau_{1}\psi=0 & \mbox{ in }\left(  0,R\right)
\\
& \\
\psi\left(  0\right)  =\psi\left(  R\right)  =0, &
\end{array}
\right.
\]
and in particular
\begin{equation}
\left\{
\begin{array}
[c]{ll}
\psi^{\prime\prime}+\psi^{\prime}\left(  \dfrac{N-1}{r}-r\right)  -\dfrac
{N-1}{r^{2}}\psi+\tau_{1}\psi\leq0 & \mbox{ in }\left(  0,r_{0}\right)  ,\\
& \\
\psi\left(  0\right)  =0,\text{ }\psi^{\prime}\left(  r_{0}\right)  \geq0. &
\end{array}
\right.  \label{prob_3}
\end{equation}
Now we multiply equation in (\ref{prob_1}) by $r^{N-1}\psi$ $\varphi_{N}$ and
equation in (\ref{prob_3}) by $r^{N-1}w_{1}$ $\varphi_{N},$ respectively.
Hence, by subtracting, we obtain%
\[
r^{N-1}\varphi_{N}(\psi w_{1}^{\prime\prime}-w_{1}\psi^{\prime\prime}
)+r^{N-1}\varphi_{N}\left(  \frac{N-1}{r}-r\right)  (\psi w_{1}^{\prime}
-w_{1}\psi^{\prime})+(\nu_{1}-\tau_{1})w_{1}r^{N-1}\psi\varphi_{N}\geq0\text{
\ \ in \ \ }(0,r_{0}).
\]
Integrating by parts the above inequality on $(0,r_{0}),$ we get%
\begin{align*}
(\nu_{1}-\tau_{1})\int_{0}^{r_{0}}w_{1}r^{N-1}\psi\varphi_{N}  &  >\int
_{0}^{r_{0}}\varphi_{N}w_{1}(r^{N-1}\psi^{\prime})^{\prime}-\varphi_{N}%
\psi(r^{N-1}w_{1}^{\prime})^{\prime}+r^{N}\varphi_{N}(\psi w_{1}^{\prime
}-w_{1}\psi^{\prime})dr\\
&  =r_{0}^{N-1}\varphi_{N}(r_{0})\left(  \psi^{\prime}(r_{0})w(r_{0}%
)-\psi(r_{0})w^{\prime}(r_{0})\right)  >0
\end{align*}
In other words
\[
\nu_{1}(R)<\tau_{1}(R)\text{ \quad}\forall R\in J_{1}.
\]

\bigskip

\textbf{Case 2}: $R\in J_{2}=\left(  \sqrt{N-1},\sqrt{N-1}+\frac{\pi}{\sqrt
{8}}\right]  .$

The above proof does not work in $J_{2}$. This because when $r>\sqrt{N-1}$ one
cannot exclude a priori that $r_{0}>\sqrt{N-1}$ too. Hence (\ref{eq_for_g''})
does no longer guarantee (\ref{g''(r_0)>0}). Clearly we may assume here that
\begin{equation}
\sqrt{N-1}<r_{0}<R, \label{Case2_r0}
\end{equation}
indeed, if not (i.e. if $r_{0}\leq\sqrt{N-1}$), we can get the claim by
repeating the arguments of Case 1.

By (\ref{prob_2}) we get
\begin{equation}
\left\{
\begin{array}
[c]{ll}
g_{1}^{\prime\prime}+\tau_{1}g_{1}<0 & \mbox{ in }(r_{0},R),\\
& \\
g_{1}\left(  r_{0}\right)  =g_{1}^{\prime}\left(  R\right)  =0. &
\end{array}
\right.  \label{g2<tau2}
\end{equation}
Multiplying the equation in (\ref{g2<tau2}) by $g_{1}\left(  r\right)  <0$ and
integrating between $r_{0}$ and $R,$ we get
\[
\int_{r_{0}}^{R}(g_{1}^{\prime})^{2}dr<\tau_{1}\int_{r_{0}}^{R}(g_{1})^{2}dr,
\]
that implies%
\[
\tau_{1}(R)>\min_{v\neq0:\text{ }v\left(  r_{0}\right)  =v^{\prime}\left(
R\right)  =0}\frac{\int_{r_{0}}^{R}(v^{\prime})^{2}dr}{\int_{r_{0}}^{R}%
v^{2}dr}=\frac{\pi^{2}}{4(R-r_{0})^{2}}.
\]
Finally, taking into account that we are under the assumption (\ref{Case2_r0}%
), we get the following
\begin{equation}
\tau_{1}(R)>\frac{\pi^{2}}{4(R-r_{0})^{2}}>\frac{\pi^{2}}{4(R-\sqrt{N-1})^{2}%
}:=h(R),\text{ \ }\forall R\in J_{2}. \label{I relazione}%
\end{equation}
Now we want to provide an estimate from above for $\nu_{1}(R),$ namely
$\nu_{1}(R)<k(R).$ To this aim we firstly note that for the values of $R\in
J_{2}$ such that $\nu_{1}(R)\leq\tau_{1}(R)$ we have%
\begin{equation}
\nu_{1}=\underset{%
\begin{array}
[c]{c}%
v\in H^{1}\left(  B_{R}\right)  ,\text{ }v\neq0\\
\int_{B_{R}}vd\gamma_{N}=0
\end{array}
}{\min}\frac{\int_{B_{R}}\left\vert Dv\right\vert ^{2}d\gamma_{N}}{\int
_{B_{R}}\left\vert v\right\vert ^{2}d\gamma_{N}}. \label{Case1}%
\end{equation}
While for the remaining values of $R$ we have to impose also the orthogonality
with $g_{1},$ that is%
\begin{equation}
\nu_{1}=\underset{%
\begin{array}
[c]{c}%
v\in H^{1}\left(  B_{R}\right)  ,\text{ }v\neq0\\
\int_{B_{R}}vd\gamma_{N}=0,\text{ }\int_{B_{R}}vg_{1}d\gamma_{N}=0
\end{array}
}{\min}\frac{\int_{B_{R}}\left\vert Dv\right\vert ^{2}d\gamma_{N}}{\int
_{B_{R}}\left\vert v\right\vert ^{2}d\gamma_{N}}. \label{Case2}%
\end{equation}
In both cases $v=x_{i}$ for $i=1,...,N$ are admissible trial functions for
$\nu_{1}$ and hence
\[
\nu_{1}\leq\frac{\gamma_{N}\left(  B_{R}\right)  }{\int_{B_{R}}x_{1}%
^{2}d\gamma_{N}},\text{ ... ,}\nu_{1}\leq\frac{\gamma_{N}\left(  B_{R}\right)
}{\int_{B_{R}}x_{N}^{2}d\gamma_{N}}.
\]
So%
\[
\frac{N}{\nu_{1}}\geq\frac{\int_{B_{R}}\left(  x_{1}^{2}+...+x_{N}^{2}\right)
d\gamma_{N}}{\gamma_{N}\left(  B_{R}\right)  };
\]
and%
\begin{equation}
\nu_{1}=\nu_{1}(R)\leq\frac{N\int_{0}^{R}e^{-\frac{s^{2}}{2}}s^{N-1}ds}%
{\int_{0}^{R}e^{-\frac{s^{2}}{2}}s^{N+1}ds}:=k\left(  R\right)  .
\label{dis v_1}%
\end{equation}
At this point we observe that $k\left(  R\right)  $ is a decreasing function,
indeed
\[
k ^{\prime}\left(  R\right)  =\frac{Ne^{-\frac{R^{2}}{2}}R^{N-1}}{\left(
\int_{0}^{R}e^{-\frac{s^{2}}{2}}s^{N+1}ds\right)  ^{2}}\left[  \int_{0}
^{R}e^{-\frac{s^{2}}{2}}s^{N+1}ds-R^{2}\int_{0}^{R}e^{-\frac{s^{2}}{2}}
s^{N-1}ds\right]  <0,
\]
where the quantity in the square brackets is negative because
\[
\int_{0}^{R}e^{-\frac{s^{2}}{2}}s^{N+1}ds=\int_{0}^{R}s^{2}e^{-\frac{s^{2}}
{2}}s^{N-1}ds<R^{2}\int_{0}^{R}e^{-\frac{s^{2}}{2}}s^{N-1}ds,
\]
Furthermore the function $h\left(  R\right)  ,$ defined in (\ref{I relazione}
), is obviously a decreasing function.

Let us consider the case $N=2$ first. Let $\overset{\_}{R}$ be the unique
positive zero of the function $f(t)=t^{2}+1-e^{\frac{t^{2}}{2}}$
\ ($\overset{\_}{R}\simeq1.585$). If $1<R<\overset{\_}{R}$ then by
(\ref{I relazione}), (\ref{dis v_1}) and by the monotonicity properties of the
functions $k(R)$ and $h\left(  R\right)  $ we get

\begin{equation}
\nu_{1}(R)<k(R)<\underset{\left(  1,\overset{\_}{R}\right)  }{\sup
}k(R)=k(1),\text{ }\forall R\in(1,\overset{\_}{R}) \label{vu1<}
\end{equation}
and%
\begin{equation}
h(\overset{\_}{R})=\underset{\left(  1,\overset{\_}{R}\right)  }{\inf
}h(R)<h\left(  R\right)  <\tau_{1}(R),\text{ }\forall R\in(1,\overset{\_}{R}).
\label{tau2>}
\end{equation}
Now since
\begin{equation}
k(1)=\frac{2\int_{0}^{1}e^{-\frac{t^{2}}{2}}tdt}{\int_{0}^{1}e^{-\frac{t^{2}
}{2}}t^{3}dt}=\frac{2-2e^{-\frac{1}{2}}}{2-3e^{-\frac{1}{2}}}\simeq4.362
\label{k(1)}
\end{equation}
and
\begin{equation}
h(\overset{\_}{R})=\frac{\pi^{2}}{4(\overset{\_}{R}-1)^{2}}\simeq7.\,210,
\label{h(R_bar)}
\end{equation}
taking into account of (\ref{vu1<}) and (\ref{tau2>}), we get
\[
\nu_{1}(R)<\tau_{1}(R),\text{ }\forall R\in(1,\overset{\_}{R}).
\]
Let us consider the remaining interval $\left[  \overset{\_}{R},1+\frac{\pi
}{\sqrt{8}}\right]  $. Since
\[
k\left(  \overset{\_}{R}\right)  =\frac{2-2e^{-\frac{\overset{\_}{R}}{2}}
}{2-3e^{-\frac{\overset{\_}{R}}{2}}}\simeq1.705<h\left(  1+\frac{\pi}{\sqrt
{8}}\right)  =2,
\]
arguing as before we get

\[
\nu_{1}(R)<\tau_{1}(R),\text{ }\forall R\in\left[  \overset{\_}{R},1+\frac
{\pi}{\sqrt{8}}\right]  .
\]

Now let $N\geq3$. If $\sqrt{N-1}<R<\sqrt{N+2},$ by (\ref{I relazione}) and
(\ref{dis v_1}) we get

\begin{equation}
\nu_{1}(R)<k(R)<\underset{\left(  \sqrt{N-1},\sqrt{N+2}\right)  }{\sup
}k(R)=k\left(  \sqrt{N-1}\right)  . \label{vu1<k}
\end{equation}
We claim that
\begin{equation}
k\left(  \sqrt{N-1}\right)  \leq\frac{2N+1}{N-1}. \label{k(N_1)}
\end{equation}
Indeed by an integration by parts the claim becomes
\[
k\left(  \sqrt{N-1}\right)  =\frac{\left(  N-1\right)  ^{\frac{N}{2}}
e^{-\frac{N-1}{2}}+\int_{0}^{\sqrt{N-1}}e^{-\frac{s^{2}}{2}}s^{N+1}ds}
{\int_{0}^{\sqrt{N-1}}e^{-\frac{s^{2}}{2}}s^{N+1}ds}\leq\frac{2N+1}{N-1}.
\]
In order to prove the above inequality it suffices to show that
\begin{equation}
e^{\frac{N-1}{2}}\int_{0}^{\sqrt{N-1}}e^{-\frac{s^{2}}{2}}s^{N+1}ds\geq
\frac{\left(  N-1\right)  ^{\frac{N}{2}+1}}{N+2}. \label{equiv_claim}
\end{equation}
Inequality (\ref{equiv_claim}), and hence the claim (\ref{k(N_1)}), easily
follows by observing that
\[
e^{\frac{N-1}{2}}\int_{0}^{\sqrt{N-1}}e^{-\frac{s^{2}}{2}}s^{N+1}
ds>e^{\frac{N-1}{2}}\int_{0}^{\sqrt{N-1}}s^{N+1}ds=e^{\frac{N-1}{2}}
\frac{(N-1)^{\frac{N+2}{2}}}{N+2}>\frac{\left(  N-1\right)  ^{\frac{N}{2}+1}
}{N+2}.
\]
Now we want to prove that
\begin{equation}
\frac{2N+1}{N-1}<h\left(  \sqrt{N+2}\right)  =\frac{\pi^{2}}{4(\sqrt
{N+2}-\sqrt{N-1})^{2}},\text{\ \ }\forall N\geq3. \label{<h}
\end{equation}
It is elementary to verify that (\ref{<h}) is true for $N=3.$ On the other
hand observe that (\ref{<h}) is false for $N=2$, that is the reason we were
forced to split Case 2 in the proof of Lemma \ref{Lemma}\ in these subcases.

Finally we get (\ref{<h}) since the sequences $\frac{2N+1}{N-1}$ and
$\frac{\pi^{2}}{4(\sqrt{N+2}-\sqrt{N-1})^{2}}$ are decreasing and increasing
respectively. Therefore, from the monotonicity of the functions $k(R)$ and
$h(R),$ (\ref{vu1<k}), (\ref{k(N_1)}) and (\ref{<h}) yield

\[
\nu_{1}(R)<k(R)<k(\sqrt{N-1})\leq\frac{2N+1}{N-1}<h\left(  \sqrt{N+2}\right)
<h\left(  R\right)  <\tau_{1}(R),\text{ }\forall R\in\left(  \sqrt{N-1}
,\sqrt{N+2}\right)  .
\]
Finally let $R\in\left[  \sqrt{N+2},\sqrt{N-1}+\frac{\pi}{\sqrt{8}}\right]  $.
We claim that
\begin{equation}
k\left(  \sqrt{N+2}\right)  \leq2. \label{k<2}
\end{equation}
Indeed arguing as before we have%
\begin{align*}
k\left(  \sqrt{N+2}\right)  -2  &  =\frac{\left(  N+2\right)  ^{\frac{N}{2}
}e^{-\frac{N+2}{2}}-\int_{0}^{\sqrt{N+2}}e^{-\frac{s^{2}}{2}}s^{N+1}ds}
{\int_{0}^{\sqrt{N+2}}e^{-\frac{s^{2}}{2}}s^{N+1}ds}\\
&  <\frac{\left(  N+2\right)  ^{\frac{N}{2}}e^{-\frac{N+2}{2}}-\int_{0}
^{\sqrt{N+2}}s^{N+1}ds}{\int_{0}^{\sqrt{N+2}}e^{-\frac{s^{2}}{2}}s^{N+1}
ds}=\frac{\left(  N+2\right)  ^{\frac{N}{2}}\left(  e^{-\frac{N+2}{2}
}-1\right)  }{\int_{0}^{\sqrt{N+2}}e^{-\frac{s^{2}}{2}}s^{N+1}ds}<0.
\end{align*}

Finally we have

\begin{align*}
\nu_{1}\left(  R\right)   &  <k(R)<k\left(  \sqrt{N+2}\right)  <2=h\left(
\sqrt{N-1}+\frac{\pi}{\sqrt{8}}\right) \\
&  <h\left(  R\right)  <\tau_{1}\left(  R\right)  ,\text{ \ }\forall
R\in\left(  \sqrt{N+2},\sqrt{N-1}+\frac{\pi}{\sqrt{8}}\right]  .
\end{align*}

\textbf{Case 3}: $R\in J_{3}=\left(  \sqrt{N-1}+\frac{\pi}{\sqrt{8}}
,+\infty\right)  .$

Before addressing this last case let us remark that the above method cannot be
used for large values of $R$. Indeed when $N=2$, for instance, we have
\[
\lim\limits_{R\rightarrow+\infty}k(R)=\lim\limits_{R\rightarrow+\infty}
\frac{2-2e^{-\frac{R^2}{2}}}{2-(R^2+2)e^{-\frac{R^2}{2}}}=1\text{ and }\lim
\limits_{R\rightarrow+\infty}h(R)=0.
\]
Therefore the inequality $k(R)<h(R),$ we have used in Case 2, does not hold
for any $R\in J_{3}$.

In order to analyze the problem for large value of the radius $R$ it appears
natural to consider the solution to problem (\ref{prob_2}) with $R=+\infty.$
Its first radial eigenfunction, as well known, is
\[
g_{\infty}\left(  r\right)  =\sum\limits_{i=1}^{N}H_{2}(x_{i})=r^{2}-N,
\]
where $H_{2}$ is the Hermite polynomial defined in (\ref{Hermite}). More
explicitly we have

\begin{equation}
\left\{
\begin{array}
[c]{ll}
g_{\infty}^{^{\prime\prime}}+g_{\infty}^{\prime}\left(  \dfrac{N-1}
{r}-r\right)  +\tau_1(\infty)g_{\infty}=0 & \mbox{ in }(0,+\infty)\\
& \\
g_{\infty}^{\prime}\left(  0\right)  =\underset{r\rightarrow+\infty}{\lim
}\left(  g_{\infty}^{\prime}\left(  r\right)  e^{-\frac{r^{2}}{2}}\right)
=0, &
\end{array}
\right.  \label{prob_infty}
\end{equation}
where $\tau_1(\infty)=2$.

Let us denote, according to the notation used in Section 3, with
$\lambda_{1}\left(B_{r}\right)$ the first  eigenvalue of the problem
\begin{equation*}
\left\{
\begin{array}
[c]{ll}
-\sum\limits_{i=1}^{N}\dfrac{\partial}{\partial x_{i}}\left(  \varphi_{N}\left(  x\right)
\dfrac{\partial u}{\partial x_{i}}\right)  =\lambda \varphi_{N}\left(  x\right)
u & \text{ in }B_{r},\\
& \\
u=0 & \text{ on }\partial B_{r}.
\end{array}
\right.
\end{equation*}
We claim that
\begin{equation}
\tau_{1}\left(  R\right)  >\tau_{1}\left(  \infty\right).\label{claim_J3}
\end{equation}
To  this aim we may assume that $r_{0}\geq\sqrt{N}$.
Indeed if  $r_{0}<\sqrt{N}$
we have
\[
\tau_{1}(R)=\lambda_{1}\left(  B_{r_{0}}\right)  >\lambda_{1}\left(
B_{\sqrt{N}}\right)  =2=\tau_{1}(\infty).
\]

Now multiplying the equation in problem (\ref{prob_2}) by $r^{N-1}$ $\varphi
_{N}g_{\infty}$ and equation in problem (\ref{prob_infty}) by $r^{N-1}$
$\varphi_{N}g_{1}$ respectively and hence subtracting, we get

\[
r^{N-1}\varphi_{N}(g_{\infty}g_{1}^{\prime\prime}-g_{\infty}^{\prime\prime
}g_{1})+r^{N-1}\varphi_{N}\left(  \frac{N-1}{r}-r\right)  (g_{\infty}
g_{1}^{\prime}-g_{1}g_{\infty}^{\prime})+(\tau_{1}-\tau_1(\infty))g_{\infty
}g_{1}r^{N-1}\varphi_{N}=0\text{\ in }(r_{0},R).
\]

Integrating between $r_{0}$ and $R$, we get

\begin{align*}
(\tau_1(\infty)-\tau_{1})\int_{r_{0}}^{R}g_{\infty}g_{1}r^{N-1}\varphi_{N}dr  &
=\int_{r_{0}}^{R}\varphi_{N}g_{\infty}(r^{N-1}g_{1}^{\prime})^{\prime}%
-\varphi_{N}g_{1}(r^{N-1}g_{\infty}^{\prime})^{\prime}+r^{N}\varphi_{N}%
(g_{1}g_{\infty}^{\prime}-g_{\infty}g_{1}^{\prime})dr\\
&  =-\varphi_{N}(r_{0})r_{0}^{N-1}\left(  r_{0}^{2}-N\right)  g_{1}^{\prime}\left(
r_{0}\right)  -2R^{N}\varphi_{N}(R)g_{1}\left(  R\right)  >0
\end{align*}
The last inequality, since we are assuming that $r_{0}\geq\sqrt{N}$, implies
the claim (\ref{claim_J3}).

Now, recalling that $k(R)$ is a decreasing function, from (\ref{k<2}), we deduce
\[
k(R)<k\left(  \sqrt{N-1}+\frac{\pi}{\sqrt{8}}\right)  <k(\sqrt{N+2}
)\leq2\text{, \ \ }\forall R>\sqrt{N-1}+\frac{\pi}{\sqrt{8}}.
\]
The last inequalities and (\ref{dis v_1}) imply

\[
\nu_{1}(R)<k(R)<2<\tau_{1}(R)\quad\text{for }R>\sqrt{N-1}+\frac{\pi}{\sqrt{8}
}.
\]

\ \ \ \ \ \ \ \ \ \ \ \ \ \ \ \ \ \ \ \ \ \ \ \ \ \ \ \ \ \ \ \ \ \ \ \ \ \ \ \ \ \ \ \ \ \ \ \ \ \ \ \ \ \ \ \ \ \ \ \ \ \ \ \ \ \ \ \ \ \ \ \ \ \ \ \ \ \ \ \ \ \ \ \ \ \ \ \ \ \ \ \ \ \ \ \ \ \ \ \ \ \ \ \ \ \ \ \ \ \ \ \ \ \ \ \ \ \ \ \ \ \ \ \ \ \ \ \ \ \ \ \ \ \ \ \ \ \ \ \ \ $\square
$

\begin{remark}
\emph{Note that the upper bound for }$\mu_{1}(R)$\emph{ given in
(\ref{dis v_1}) is asymptotically sharp, as }$R$\emph{ goes to }$+\infty
.$\emph{ Indeed, as it is easy to verify, it holds}%
\[
\lim_{R\rightarrow+\infty}\nu_{1}(R)=\frac{N\int_{0}^{+\infty}e^{-\frac{s^{2}%
}{2}}s^{N-1}ds}{\int_{0}^{+\infty}e^{-\frac{s^{2}}{2}}s^{N+1}ds}=1=\mu
_{1}(\mathbb{R}^{N}),\text{ \ }\forall N\in\mathbb{N}.
\]

\end{remark}

\bigskip

Lemma \ref{Lemma} ensures that the first eigenfunction associated to the first
eigenvalue of problem (\ref{P2_I}) with $\Omega=B_{R},$ is in the form
$u(x)=w(\left\vert x\right\vert )Y(\theta),$ where $\theta$ belongs to
$\mathbb{S}_{1}^{N-1}$ and the radial function $w$ has one sign in $B_{R}$ and
it satisfies the following problem%
\begin{equation}
\left\{
\begin{array}
[c]{lll}%
w^{\prime\prime}(r)+w^{\prime}(r)\left(  \dfrac{N-1}{r}-r\right)  +\mu
_{1}(B_{R})w(r)-\dfrac{N-1}{r^{2}}w(r)=0, & \text{for } & r\in(0,R)\\
&  & \\
w\left(  0\right)  =w^{\prime}\left(  R\right)  =0. &  &
\end{array}
\right.  \label{eq_BR}%
\end{equation}
Multiplying the equation in ($\ref{eq_BR})$ by $w$ $\varphi_{N}$ and
integrating over $B_{R,},$ we get%
\begin{align*}
&  \mu_{1}(B_{R})\int_{B_{R}}w(\left\vert x\right\vert )^{2}d\gamma_{N}\\
&  =-N\omega_{N}\int_{0}^{R}(w^{\prime}r^{N-1})^{\prime}w(r)e^{-\frac{r^{2}
}{2}}dr+N\omega_{N}\int_{0}^{R}r^{N}w(r)w^{\prime}(r)e^{-\frac{r^{2}}{2}%
}dr+\int_{B_{R}}\frac{1}{\left\vert x\right\vert ^{2}}w\left(  \left\vert
x\right\vert \right)  ^{2}d\gamma_{N}\\
&  =\int_{B_{R}}\left(  w^{\prime}\left(  \left\vert x\right\vert \right)
\right)  ^{2}d\gamma_{N}+\int_{B_{R}}\frac{1}{\left\vert x\right\vert ^{2}
}w\left(  \left\vert x\right\vert \right)  ^{2}d\gamma_{N}.
\end{align*}
Thus%
\begin{equation}
\mu_{1}(B_{R})=\frac{\displaystyle\int_{B_{R}}\left(  \left(  w^{\prime
}\left(  \left\vert x\right\vert \right)  \right)  ^{2}+\dfrac{N-1}{\left\vert
x\right\vert ^{2}}w\left(  \left\vert x\right\vert \right)  ^{2}\right)
d\gamma_{N}}{\displaystyle\int_{B_{R}}w(\left\vert x\right\vert )^{2}
d\gamma_{N}}. \label{eig_BR}
\end{equation}
Now we are able to prove our main result.

\begin{theorem}
\label{Main_th}The ball maximizes the first Neumann eigenvalue among all
 Lipschitz open sets $\Omega$ of $\mathbb{R}^{N}$ of
prescribed Gaussian measure and symmetric about the origin. Moreover, it is the unique maximizer in this class.
\end{theorem}

\textbf{Proof} Let $B_{R}$ the ball centered at the origin having the same
Gaussian measure as $\Omega$. We define
\begin{equation}
G(r)=\left\{
\begin{array}
[c]{ll}%
w(r) & \mbox{ for }0<r<R\\
w(R) & \mbox{ for }r\geq R,
\end{array}
\right.  \label{def G}
\end{equation}
where $w$ is the solution of (\ref{prob_2}) satisfying (\ref{eig_BR}). By the
results stated above the function $G$ is nondecreasing and nonnegative. We
introduce the functions
\[
P_{i}(x)=G(\left\vert x\right\vert )\frac{x_{i}}{\left\vert x\right\vert
}\text{ \ \ for \ \ }1\leq i\leq N.
\]
The assumption on the symmetry of $\Omega$  guarantees
\begin{equation}
\int_{\Omega}P_{i}(x)d\gamma_{N}=0,\text{ \ \ }\forall
i=1,...,N.\label{Orth_P_i}
\end{equation}
Hence each function $P_{i}$ is admissible in the variational formulation
(\ref{def_auto}).

Since
\[
\frac{\partial P_{i}}{\partial x_{j}}=G^{\prime}(\left\vert x\right\vert
)\frac{x_{i}x_{j}}{\left\vert x\right\vert ^{2}}-G(\left\vert x\right\vert
)\frac{x_{i}x_{j}}{\left\vert x\right\vert ^{3}}+\delta_{ij}\frac{G(\left\vert
x\right\vert )}{\left\vert x\right\vert },
\]
where $\delta_{ij}$ is the Kronecker symbol, summing over $j=1,...,N$, we get
\begin{equation}
\mu_{1}(\Omega)\leq\frac{\displaystyle\int_{\Omega}\left(  \left(  G^{\prime
}\left(  \left\vert x\right\vert \right)  \right)  ^{2}\frac{x_{i}^{2}
}{\left\vert x\right\vert ^{2}}-G^{2}\left(  \left\vert x\right\vert \right)
\frac{x_{i}^{2}}{\left\vert x\right\vert ^{4}}+\frac{G^{2}\left(  \left\vert
x\right\vert \right)  }{\left\vert x\right\vert ^{2}}\right)  d\gamma_{N}
}{\displaystyle\int_{\Omega}G(\left\vert x\right\vert )^{2}\frac{x_{i}^{2}
}{\left\vert x\right\vert ^{2}}d\gamma_{N}}. \label{mu_1_i}
\end{equation}
Set
\[
N(r)=\left(  G^{\prime}\left(  r\right)  \right)  ^{2}+\dfrac{N-1}{r^{2}}
G^{2}\left(  r\right)
\]
and
\[
D(r)=G^{2}\left(  r\right)  .
\]
Summing up inequalities (\ref{mu_1_i}) over $i=1,...,N$, the angular dependence
drops out and we finally get
\begin{equation}
\mu_{1}(\Omega)\leq\frac{\displaystyle\int_{\Omega}\left(  \left(  G^{\prime
}\left(  \left\vert x\right\vert \right)  \right)  ^{2}+\dfrac{N-1}{\left\vert
x\right\vert ^{2}}G^{2}\left(  \left\vert x\right\vert \right)  \right)
d\gamma_{N}}{\displaystyle\int_{\Omega}G(\left\vert x\right\vert )^{2}
d\gamma_{N}}=\frac{\displaystyle\int_{\Omega}N(\left\vert x\right\vert
)d\gamma_{N}}{\displaystyle\int_{\Omega}D(\left\vert x\right\vert )d\gamma
_{N}}. \label{autovalore}
\end{equation}
It is straightforward to verify that
\[
\frac{d}{dr}N(r)<0.
\]
Now we claim that%
\begin{equation}
\int_{\Omega}N(\left\vert x\right\vert )d\gamma_{N}\leq\int_{B_{R}
}N(\left\vert x\right\vert )d\gamma_{N}. \label{dis_BR_1}
\end{equation}
Hardy-Littlewood inequality (\ref{HL}) ensures
\begin{equation}
\int_{\Omega}N(\left\vert x\right\vert )d\gamma_{N}\leq\int_{0}^{\gamma
_{N}(\Omega)}N^{\ast}(s)ds=\int_{0}^{\gamma_{N}(B_{R})}N^{\ast}(s)ds,
\label{dis_1}
\end{equation}
where $N^{\ast}$ is the decreasing rearrangement of $N$. Setting $s=\gamma
_{N}(B_{r})=\frac{N\omega_{N}}{\left(  2\pi\right)  ^{N/2}}\displaystyle\int
_{0}^{r}e^{-\frac{s^{2}}{2}}s^{N-1}ds,$ we get
\[
\int_{0}^{\gamma_{N}(B_{R})}N^{\ast}(s)ds=\frac{N\omega_{N}}{\left(
2\pi\right)  ^{N/2}}\int_{0}^{R}N^{\ast}(\gamma_{N}(B_{r}))r^{N-1}
e^{-\frac{r^{2}}{2}}dr.
\]
Note that%
\[
N^{\ast}(\gamma_{N}(B_{r}))=N(r),
\]
since $N^{\ast}(\gamma_{N}(B_{r}))$ and $N(r)$ are equimeasurable and both
radially decreasing functions. Therefore
\begin{equation}
\frac{N\omega_{N}}{\left(  2\pi\right)  ^{N/2}}\int_{0}^{R}N^{\ast}(\gamma
_{N}(B_{r}))r^{N-1}e^{-\frac{r^{2}}{2}}dr=\frac{N\omega_{N}}{\left(
2\pi\right)  ^{N/2}}\int_{0}^{R}N(r)r^{N-1}e^{-\frac{r^{2}}{2}}dr=\int_{B_{R}
}N(\left\vert x\right\vert )d\gamma_{N} \label{dis_2}%
\end{equation}
Combining (\ref{dis_1}) and (\ref{dis_2}), we obtain the claim (\ref{dis_BR_1}
). Analogously it is possible to prove that%
\begin{equation}
\int_{\Omega}D(\left\vert x\right\vert )d\gamma_{N}\geq\int_{B_{R}
}D(\left\vert x\right\vert )d\gamma_{N}. \label{dis_BR_2}
\end{equation}
Indeed since $D$ is an increasing function, we have%
\begin{align*}
\int_{\Omega}D(\left\vert x\right\vert )\varphi_{N}(\left\vert x\right\vert
)dx  &  \geq\int_{0}^{\gamma_{N}(B_{R})}D_{\ast}(s)ds\\
&  =\frac{N\omega_{N}}{\left(  2\pi\right)  ^{N/2}}\int_{0}^{R}D_{\ast
}(1-e^{-\frac{r^{2}}{2}})r^{N-1}e^{-\frac{r^{2}}{2}}dr=\int_{B_{R}
}D(\left\vert x\right\vert )\varphi_{N}(\left\vert x\right\vert )dx,
\end{align*}
where $D_{\ast}$ is the increasing rearrangement of $D$. By (\ref{def G}
),(\ref{dis_BR_1}) and (\ref{dis_BR_2}), the equality (\ref{autovalore})
becomes%
\[
\mu_{1}(\Omega)\leq\frac{\displaystyle\int_{B_{R}}\left(  \left(  w^{\prime
}\left(  \left\vert x\right\vert \right)  \right)  ^{2}+\dfrac{N-1}{\left\vert
x\right\vert ^{2}}w\left(  \left\vert x\right\vert \right)  ^{2}\right)
d\gamma_{N}}{\displaystyle\int_{B_{R}}w(\left\vert x\right\vert )^{2}%
d\gamma_{N}}=\mu_{1}(B_{R}),
\]
which is the desired inequality. Moreover, from the monotonicity properties of
the functions $N$ and $D$, it easy to realize that inequalities
(\ref{dis_BR_1}) and (\ref{dis_BR_2}) reduce to equalities only when $\Omega$
is the ball $B_{R}.$
\ \ \ \ \ \ \ \ \  \ \ \ \ \ \ \ \ \ \ \ \ \  \ \ \ \ \ \ \ \ \ \ \ \  \ \ \ \ \ \ \ \ \ \ \ \ \ \ \ \ \ \
 \ \ \ \ \ \ \ \ \    \ \ \ \ \ \ \ \ \    \ \ \ \ \ \ \ \ \   \ \ \ \ \ \ \ \ \   \ \ \ \ \ \ \ \ \  \ \ \ \ \ \ \ \ \
  \ \ \ \ \ \ \ \ \ $\square
$

\begin{remark}
Note that the assumption on the symmetry of $\Omega$ is used solely to guarantee the orthogonality conditions \emph{(\ref{Orth_P_i})}.
\end{remark}

\begin{remark}
\label{remark}Since the half-spaces are not symmetric about the origin, Theorem
\emph{\ref{Main_th}} cannot exclude the possibility that such a domain maximizes
$\mu_{1}\left(  \Omega\right)  $ in dimension greater than one. This
phenomenon does not occur since any half-space has first Neumann eigenvalue
equal to $1$, independently of its measure. It is easy to show an example of a
set which is not symmetric about the origin whose first Neumann eigenvalue is bigger than $1$.
Consider, for instance, in $\mathbb{R}^{2}$ the square $T=\left(
\sqrt{3-\sqrt{6}},\sqrt{3+\sqrt{6}}\right)  ^{2}.$ As it is immediate to
verify, $\mu_{1}\left(  T\right)  =5$ \ and it is a double eigenvalue. A
corresponding eigenfunction is $u_{1}\left(  x,y\right)  =u_{1}(x)=H_{5}
(x)=x^{5}-10x^{3}+15x.$ This simply follows by observing that $H_{5}^{\prime
}(x)<0$ $\ \forall x\in\left(  \sqrt{3-\sqrt{6}},\sqrt{3+\sqrt{6}}\right)  $
and $H_{5}^{\prime}\left(  \sqrt{3-\sqrt{6}}\right)  =H_{5}^{\prime}\left(
\sqrt{3+\sqrt{6}}\right)  =0$. Let us round a corner of this square by
considering the family of domains%
\[
T_{\delta}=\left\{  (x,y)\in\mathbb{R}^{2}:\sqrt{3-\sqrt{6}}\leq x\leq
\sqrt{3+\sqrt{6}}\text{ and }\sqrt{3-\sqrt{6}}\leq y\leq f_{\delta
}(x)\right\}  \text{,}
\]
with $\delta<<1$ and
\[
f_{\delta}(x)=\left\{
\begin{array}
[c]{ccc}
\sqrt{3+\sqrt{6}} & \text{if} & \sqrt{3-\sqrt{6}}\leq x\leq\sqrt{3+\sqrt{6}
}-\delta\\
&  & \\
\sqrt{3+\sqrt{6}}-\delta+\sqrt{\delta^{2}-\left(  x-\left(  \sqrt{3+\sqrt{6}
}-\delta\right)  \right)  ^{2}} & \text{if} & \sqrt{3+\sqrt{6}}-\delta
<x\leq\sqrt{3+\sqrt{6}}.
\end{array}
\right.
\]
Now the first non trivial Neumann eigenfunction relative to $T_{\delta}$
cannot depend on one variable only. The sequence of compact sets $T_{\delta}$
converges, in the Hausdorff distance, to $T,$ and therefore, see \cite{Che},
we have that $\mu_{1}\left(  T_{\delta}\right)  \rightarrow\mu_{1}\left(
T\right)  $. Therefore for $\delta$ small enough we have%
\[
5+O(1)=\mu_{1}(T_{\delta})>1=\mu_{1}\left(  T^{\bigstar}\right)  ,
\]
where $T^{\bigstar}$ is the half-space having the same Gaussian measure as
$T.$\bigskip
\end{remark}

\section{Appendix}

Here we want to show that $\tau_{k}(R)$, the
nontrivial eigenvalues of (\ref{prob_2}), are all decreasing functions. To
this aim we apply, in this simple case, the shape derivative formula for
Neumann eigenvalues, see, e.g., \cite{H1} and \cite{H2}.

Let $R\left(  t\right)  =R+t,$ with $t>0,$ and let$\ \mu_{k}\left(  t\right)
=\mu_{k}\left(  0,R\left(  t\right)  \right)  $ be the $k$-th eigenvalue of
problem
\begin{equation}
\left\{
\begin{array}
[c]{lll}%
-u_{rr}\left(  r,t\right)  +ru_{r}\left(  r,t\right)  -\dfrac{N-1}{r}%
u_{r}\left(  r,t\right)  =\mu_{k}\left(  t\right)  u\left(  r,t\right)  &
\text{in} & (0,R\left(  t\right)  )\\
&  & \\
\left.  u_{r}\left(  r,t\right)  \right\vert _{r=0,R\left(  t\right)  }=0, &
&
\end{array}
\right.  \label{prob_Appendix1}
\end{equation}
and, finally, let $u\left(  r,t\right)  $ be an eigenfunction corresponding to
$\mu_{k}\left(  t\right)  $ such that%
\begin{equation}
\left\Vert u\right\Vert _{L^{2}\left(  B_{R}\left(  t\right)  ,\gamma
_{N}\right)  }^{2}=\frac{N\omega_{N}}{\left(  2\pi\right)  ^{N/2}}\int
_{0}^{R\left(  t\right)  }u^{2}\left(  r,t\right)  r^{N-1}e^{-\frac{r^{2}}{2}%
}dr=1. \label{cond_normal}
\end{equation}
We have

\begin{proposition}
\label{Deriv_atuoAppendix}It holds that
\begin{equation}
\mu_{k}^{\prime}\left(  0\right)  =-\frac{N\omega_{N}}{\left(  2\pi\right)
^{N/2}}\mu_{k}\left(  0\right)  u^{2}\left(  R\right)  R^{N-1}e^{-\frac{R^{2}
}{2}} \label{Derivata_auto_}
\end{equation}
where $u\left(  r\right)  =u\left(  r,0\right)  $ is the eigenfunction of
problem (\ref{prob_Appendix1}) in $(0,R).$
\end{proposition}

\textbf{Proof. }Differentiating (\ref{cond_normal}) we have for $t=0,$
\begin{equation}
2\int_{0}^{R}u\left(  r\right)  u_{t}\left(  r,0\right)  r^{N-1}
e^{-\frac{r^{2}}{2}}dr=-e^{-\frac{R^{2}}{2}}R^{N-1}u^{2}\left(  R\right)  .
\label{XXX}
\end{equation}

Multiplying the equation in (\ref{prob_Appendix1}) by $u\left(  r,t\right)  $
$r^{N-1}e^{-\frac{r^{2}}{2}},$ we get
\[
\mu_{k}\left(  t\right)  u^{2}\left(  r,t\right)  r^{N-1}e^{-\frac{r^{2}}{2}
}=-(r^{N-1}u_{r}\left(  r,t\right)  )_{r}\text{ \ }u\left(  r,t\right)
e^{-\frac{r^{2}}{2}}+r^{N}u\left(  r,t\right)  u_{r}\left(  r,t\right)
e^{-\frac{r^{2}}{2}}.
\]
Integrating the above equality on $(0,R\left(  t\right)  )$ and recalling
condition (\ref{cond_normal}) we get
\begin{equation}
\mu_{k}\left(  t\right)  =\frac{N\omega_{N}}{\left(  2\pi\right)  ^{N/2}}
\int_{0}^{R\left(  t\right)  }u_{r}^{2}\left(  r,t\right)  r^{N-1}
e^{-\frac{r^{2}}{2}}dr. \label{eq_Appe}
\end{equation}
Differentiating we obtain
\[
\mu_{k}^{\prime}\left(  0\right)  =\frac{2N\omega_{N}}{\left(  2\pi\right)
^{N/2}}\int_{0}^{R}u_{r}\left(  r,0\right)  u_{rt}\left(  r,0\right)
r^{N-1}e^{-\frac{r^{2}}{2}}dr=\frac{2N\omega_{N}}{\left(  2\pi\right)  ^{N/2}
}\mu_{k}\left(  0\right)  \int_{0}^{R}u\left(  r\right)  u_{t}\left(
r,0\right)  r^{N-1}e^{-\frac{r^{2}}{2}}dr.
\]
So by (\ref{XXX}) we obtain the claim (\ref{Derivata_auto_}%
).\ \ \ \ \ \ \ \ \ \ \ \ \ \ \ \ \ \ \ \ \ \ \ \ \ \ \ \ \ \ \ \ \ \ \ \ \ \ \ \ \ \ \ \ \ \ \ \ \ \ \ \ \ \ \ \ \ \ \ \ \ \ \ \ \ \ \ \ \ \ \ \ \ \ \ \ \ \ $\square
$

\bigskip

\textbf{Acknowledgments.} The authors would like to thank Prof. Mark S.
Ashbaugh for helpful comments and suggestions.

\end{document}